\documentclass[14pt]{extarticle}
\usepackage[margin=1in]{geometry}
\usepackage{dsfont}
\usepackage{amsfonts}
\usepackage{tikz}
\usetikzlibrary{matrix}

\newtheorem{theorem}{Theorem}[section]
\newtheorem{pro}{Proposition}[section]
\newtheorem{lemma}{Lemma}[section]

\newtheorem{cor}{Corollary}[section]

\newcommand{\proof}[1]{\noindent{\it\bf Proof:#1\ }}
\newcommand{\QED}{\hfill$\Box$\medskip}

\begin{document}

\title{   Higher-degree Smoothness of Perturbations  II }
\author{Gang Liu\\Department of Mathematics\\UCLA  }
\date{September  1,  2016}
\maketitle

 \section{  Introduction}
 In this paper we  generalize part of the higher-degree  smoothness  results  in  perturbation theory in \cite {3} from the case that  the   stable maps   have the  fixed domain $S^2$ to the general genus zero  case. Note that   genus zero case already captures all   analytic difficulties related to  the lack of differentiability of transition functions between local slices (see \cite {6} for  the  discussion on this ). 
 The results in this paper and  \cite {4} together give  one of the two  methods for the infinite dimensional  set-up used in \cite {5}
 (compare the other infinite dimensional  set-up in \cite{1}).

 
 The main  result  of this paper  is  the following theorem (see the relevant definitions in the later sections).
 
 \begin{theorem}
 	Let $K\simeq K_t$ with $t\in { \bar  W}(\Sigma)$ be the fixed part of the local universal family of stable curves ${\cal S}\rightarrow N(\Sigma)$, where 
 	$N(\Sigma)$ is a small open  neighborhood of $[{\Sigma}]$  in  ${\overline {\cal M}}_{0, k}$ with the  local coordinate chart  ${\bar W}(\Sigma)$. Consider   the  local uniformizer (slice)
 	 ${ W}(f, {\bf H}_f)$ 
 	 centered at $f:\Sigma\rightarrow M$ of stable $L_k^p$-maps with domains ${\cal S}_t, t\in {\bar W}(\Sigma)$ and  the corresponding space ${\widetilde W}(f_{K})$ of $L_k^p$ maps with domain $K$ with associated bundle ${\cal L}^{K}\rightarrow {\widetilde W}(f_{K})$. 
 	Let ${\xi}^{K}:{\widetilde W}(f_{K})\rightarrow {\cal L}^{K}$ be a smooth section satisfying the condition $C_1$ and $C_2$. Then  ${\xi}^{K}$   gives rise   a stratified  $C^{m_0}$-smooth section   $\xi: { W}(f, {\bf H}_f)\rightarrow {\cal L}$  on the local slice  ${ W}(f, {\bf H}_f)$, which,  viewed  in any other  local slice, 
 	is still  stratified $C^{m_0}$-smooth   on their "common intersections
    " (=the fiber product over the space of unparametrized stable maps). 
 	
 \end{theorem}
 
 Here  the conditions $C_1$ and $C_2$ are defined in \cite {3} using the bi-grading there 
 (see the definition in  \cite {3} ) as follows.

 $C_1:$ The section   $\eta:S_f\rightarrow  {\cal L}|_{S_f}$ can be extended into  a $C^{m_0}$-smooth section
 $\eta_{-m}:(S_f)_{-m}\rightarrow  {\cal L}_0$ for some $m\geq m_0$
 
 $C_2:$ The image of  $\eta(h)$ is  lying in $L_{k+m}^p(\Sigma, E_h)=:({\cal L}_h)_{m}$ with $m\geq m_0$ so that $\eta_{-m, m}:(S_f)_{-m}\rightarrow  {\cal L}_m$ is smooth.
 
 Recall that $m_0=[k-2/p]$.  We will assume that $p>2$ and $m_0>1$ throughout this paper as in   \cite {3}.

 This theorem is  proved in Sec.5. The two kinds of Banach  neighborhoods on  an  end
 near a stable nodal map are defined in section 3.  The corresponding (stratified ) smooth structures on each of such neighborhoods  are defined in section 4.  

Only elementary facts on Sobolev spaces and standard calculus on Banach spaces are used in this paper, for which  we refer to \cite{2, 7}.



 	

 \section {Local universal family of stable curves}
 
The starting point of this paper is the local deformation theory of stable maps.  To this end, we 
 need  recall  the local deformation of the stable curves first.

\subsection{Stable curves}
 
 Given an "initial" stable curve $\Sigma^0$, let $T_0$ be the tree associated to the domain  of the stable curve   $\Sigma^0=(S^0, {\bf d}^0,  {\bf x}^0)$ so that the underlying curve $S^0$ is a nodal surface with desingularization ${\hat S^0}=\coprod_{v\in T_0}S^0_v$ as the disjoint union of its components labeled by the vertices $v\in T_0$.
Here the  double points ${\bf d}^0=\cup_{v\in T_0}{\bf d}^0_v$ with each ${\bf d}^0_v=\{ d^0_{vw}, [vw]\in E(T_0)\}$, where each double point  $ d^0_{vw}$ on  $S^0_v$ is labeled by an edge in  the set of edges $E(T_0)$ of $T_0$; the marked points ${\bf x}^0=\cup_{v\in T_0}{\bf x}^0_v$ with  ${\bf x}^0_v$  to be the marked points on $S^0_v$. Clearly the nodal surface $S^0$ is obtained form ${\hat S^0}$ by identifying the double points.

 In above, we have abused notations using ${\bf d}$   to denote both double points and the  set of their collections. Similarly for ${\bf x}$, and we will continue do so for other similar notations.

The  distinguished points (=the double points and marked points) of $\Sigma^0$  on $S^0_v$/$S^0$, will be denoted by ${\bf p}^0_v$/${\bf p}^0$.  Note that  the stable  curve $\Sigma^0$ determines and is determined by 
${\bf p}^0$ upto the  actions of $G(=\Pi_{v\in T_0}SL(2, {\bf C})_v)$.


 Let $N^{T_0}(\Sigma^0)$  be a small neighborhood of $[\Sigma^0]$ in the moduli space ${\cal M}^{T_0}_{0, k}$ with fixed topological type given by $T_0$, where $[\Sigma^0]\in {\cal M}^{T_0}_{0, k}$ is the "moduli point" of $\Sigma^0$. Then $N^{T_0}(\Sigma^0)$  parametrizes 
 the stable curves near $\Sigma_0$, or equivalently the nearby  distinguished points ${\bf p}$ on the same fixed ${\hat S^0}$. Hence we may introduce the parameters
  $b=\{b_{vu}; v\in T_0, p^0_{vu}\in {\bf p}^0_v\}$ with $b_{vu}$ in a small disc $D_{vu}({ p}^0_{vu})$ on $S^0_v$ centered at ${ p}^0_{vu}$.  In order to   quotient out  the (local) actions of $SL(2, {\bf C})_v$, for each $v\in T_0$, we fix the last three parameters in $b_{vu}$. Note  this also selects the three corresponding distinguished pints on $S_v$ that makes  it  {\bf marked} so that  the identification $S_v\simeq S^2$ is specified.

  Then the parameter
  $b$ is corresponding  to the stable curve $\Sigma_b=(S_b, {\bf p}_b)$. The collection of such parameters $b$ will be denoted by $W^{T_0}(\Sigma^0)$, considered as one of the natural holomorphic  coordinate charts of  $N^{T_0}(\Sigma^0)$. Of course, different choices of  fixing  three elements in each set $b_v, v\in T_0$ above give other but 
  same kind of coordinate charts of $N^{T_0}(\Sigma^0)$.

  In this notation, the initial surface, $\Sigma^0=\Sigma_0$  or  $\Sigma_{b}  $ with $b=0.$   In the following, the notations $\Sigma^0$ are used interchangeably with  $ \Sigma_0$ so that  $N^{T_0}(\Sigma^0)$, $W^{T_0}(\Sigma^0)$ etc. will also be denoted by 
  $N^{T_0}(\Sigma_0)$, $W^{T_0}(\Sigma_0)$ accordingly. Similarly   
  ${\bf p}^0_v$/${\bf p}^0={\bf p}_{0;v}/{\bf p}_0$.
  
  Note that  the desingularization ${{\hat S_b}}$ is same as  ${{\hat S^0}}=
  {\hat S_0}$ so that
  $S_b$ has the same components as $S=S_0$ has. As before,   $ {\bf p}_b=\cup_{v\in T_0}{\bf p}_{b;v}$, and  ${\bf p}_{b;v}={\bf x}_{b;v}\cup
  {\bf d}_{b;v}$ lying on $S^0_v$.
  

 Let   ${\bar N}^{T_1}(\Sigma_0)$  be a small (full) neighborhood of $\Sigma_0$ in the moduli space ${\overline {\cal M}}^{T_1}_{0, k}$, and ${ N}^{T_1}(\Sigma_0)$ be its top stratum as an open set in stratum ${ {\cal M}}^{T_1}_{0, k}$ with fixed topological type given by $T_1$. Then ${\bar N}^{T_1}(\Sigma_0)$  parametrizes 
 the stable curves $\Sigma_t$ near $\Sigma_0$ whose topological types are "bewteen $T_0$ and $T_1$".  In particular,
 when $T_1$ is the top stratum of ${\cal M}_{0, k}$, ${\bar N}^{T_1}(\Sigma_0)$  parametrizes 
 the stable curves $\Sigma_t$ of all types near $\Sigma_0. $ Here $t\in 
 {\bar W}^{T_1}(\Sigma_0)$ where $ 
 {\bar W}^{T_1}(\Sigma_0)$ is one of the natural coordinate  charts of ${\bar N}^{T_1}(\Sigma_0)$ extending $ 
 { W}^{T_0}(\Sigma_0)$.
 Thus each parameter $t\in {\bar N}^{T_1}(\Sigma_0)$ has the form $t=(b, a)$ with $b\in { W}^{T_0}(\Sigma_0)$.  Here  $a=\{a_{vw}; [vw]\in E(T_0)\}$ is the collection of the gluing parameters describes the gluing pattern from $S_b$ to the glued surface $S_t$ defined below. The non-zero entries of $a$, denoted by 
 $a^e=\{a_{vw}; [vw]\in E(T_0), C(v, w)=u\in T_{1}\}$ are the effective parameters  $a_{vw}$ gluing the components $S_v$ and $S_w$. 
 Here
  the map $C:T_0\times T_0\rightarrow T_{1}$ is  partially defined on $ T_0\times T_0$, and   for each $(v, w)$ with 
  $(v, w)\in E(T_0)$  it  is defined by  $C(v, w)=u\in T_{1}$ if  the component
  $S_{t; u}$ of $S_t$ is obtained  from the components $S_{b;v}$ and $S_{b;w}$.

   Recall the definition of $\Sigma_t=(S_t; {\bf p}_t)$ with ${\bf p}_t={\bf d}_{t}\cup {\bf x}_{t}$ as follows.
   
   (1) $S_{t}$ is obtained from $S_b$ by gluing at those double points $d_{vw}$ with 
   $a_{vw}\not=0$. Hence a component 
   $S_{t; u}=\#_{\{a_u=\{a_{v_i, v_j}\};C(v_i, v_j)=u\in T_1\}} (S_{v_1}, \cdots S_{v_{k(u)}})$. Here the  right-hand side above is the gluing of the components 
 $S_{v_1}, \cdots S_{v_{k(u)}}$ in $\Sigma_b$ with  the gluing parameter $a_u$. 
 
  For each $a_{v_i, v_j}\not =0,$ denote $(v_i, v_j)$ by $(u_+, v_{-})$ and 
 $a_{v_i, v_j}$ by $a$ temporarily. Let $D_{\pm}$ be the small discs on $S_{\pm}\simeq S^2$ with complex coordinate $w_{\pm}$. Then the gluing $\#_{a}(D_{-}, D_+)=D_{-}\coprod D_+$ quotient out the relation that $w_{-}\cdot w_+=a.$ Applying this to each nonzero $a_{v_i, v_j}$ above  gives the desired gluing.
 
  (2) The double points on 
 $\Sigma_t=\Sigma_{b,a}$ are exactly the part of the double points $ {\bf d}_b=\{{\bf d }_{b;uu'}\}$ such that 
 $a_{uu'}=0.$  Since we only consider the local deformations, we may assume that
 $|a|$ is sufficient small so that the marked points ${\bf x}_b$ become the corresponding ones, denoted by  ${\bf x}_t$ through  the gluing.
 
 \medskip
 \noindent
 ${\bullet }$   ${\bullet }$  ${\bullet }$ Fixed part $K_t\simeq K_0$ in $S_t$.
 
The  "fixed part"  $K_{\epsilon, t,u}$ of $S_{t;u}$ defined by  $$K_{\epsilon,t;u}=S_u\setminus \{\cup_{a_{v_i, v_j}\not = 0}N_{\epsilon}(b, a_{v_i, v_j})\cup_{a_{v_i, v_j}= 0}D_{\epsilon}( d_{b:v_i, v_j})\}.$$ Here $N_{\epsilon}(b, a_{v_i, v_j})$ is the "neck" part near the double point $d_{b;v_i, v_j}$ obtained by gluing the two corresponding small discs $D_{\epsilon}( d_{b;v_i, v_j})$ and $D_{\epsilon}( d_{b;v_j, v_i})$ with gluing parameter
 $a_{v_i, v_j}\not = 0$. Thus $K_{\epsilon, t,u} \subset S_{t;u}$  becomes a {\bf fixed } subset of $S_{b;u}$ independent  of $a$ with $t=(b, a)$.

 To get rid of the $b$-dependency of the fixed part,  go back to  the desingularization ${\hat S}$ of initial underlying curve $S(=S^0)$ of $\Sigma^0(=\Sigma_0)$. For each double point $d^0_{vu}$ on the component $S_v$, choose another small disc $D_{\epsilon_2}(d_{vu})$ of radius  $\epsilon_2>\epsilon$ such that $D_{\epsilon}(d_{vu}(b))\subset D_{\epsilon_2}(d_{vu})$ for all $b\in W^{T_0}(\Sigma_0)$. Then define the ( "smallest") fixed part $K$ to be the complement of the union of all the discs $D_{\epsilon_2}(d_{vu})$ in $S$. Then first of all, $K$ can be considered  as a subset in $S_b$, denote by $K_b$,  since both $S$ and $S_b$ have the same components; secondly by the construction of the gluing, it  can also be considered as a subset of $S_t$, denoted by $K_t$. Note that for $|a|$ small enough, the marked points ${\bf x}_t$ are lying on   $K_{t}$.

 Thus  for $|a|$ small enough with the types between $T_0$ and  $T_1$, $K_{t}\simeq K_{b}
\simeq K_{0}\subset S$ as the fixed part independent of $t$ while  $\{S_{t}\}$ is a family of curves  that are deforming.

 \subsection{Local universal family of the first kind}
 It is well-known that the total family obtained from the gluing construction above, ${\cal S}=:{\cal S}(\Sigma_0)\rightarrow {\bar N}^{T_1}(\Sigma_0)$ with the fiber ${\cal S}_t=S_{t}$ is a proper morphism of complex    manifolds/orbifords. 

  \begin{lemma}
  	Given $T_1\geq T_0$, let $t_0=(b_0, a_0)$ be the center of  $W^{T_1}(\Sigma_{t_0}).$ Then 	there is a smooth but non-holomorphic family of  identifications 
  	$\{\lambda^{t_0}_t:(S_t, {\bf x}_t)\rightarrow (S_{t_0}, {\bf x}_{t_0})\}$ for  $t\in W^{T_1}(\Sigma_{t_0})$.  Under this identification, the smallest fixed part 
  	$K_t$ is identified  with $K_{t_0}$ so that  $K_{t}\simeq K_{t_0}\simeq K_b\simeq K_{b_0}\simeq K_0$, the  small disks  or  "neck" areas at or near double points on $S_t$ are identifies with the correspoding ones on $S_{t_0}$. Away from the small annuli of the tubular neighborhoods of the boundaries of $K_t$,  $\lambda^{t_0}_t$ is holomorphic and preserves any of the natural metrics.

  	
  	Moreover, these maps  together give rise a smooth map
  	$\lambda^{t_0}:{\cal S}|_{N^{T_1}(\Sigma_{t_0})}\rightarrow {\cal S}_{t_0}=S_{t_0}$ and hence the induced smooth  the map ${\hat \lambda^{t_0}}=( \lambda^{t_0}, \pi):{\cal S}|_{N^{T_1}(\Sigma_{t_0})}\rightarrow {\cal S}_{t_0}\times N^{T_1}(\Sigma_{t_0})$. Here  $\pi:{\cal S}|_{N^{T_1}(\Sigma_{t_0})}\rightarrow N^{T_1}(\Sigma_{t_0})$ is the projection map.

  \end{lemma}
  
 The proof of this lemma is the immediate consequence of the construction of these diffeomorphisms below.
  
  The map 	$\lambda^{t_0}$ above will be used to define  the  smooth structure of the {\bf first kind} on the corresponding neighborhood of the {\bf first kind} of  a  stable map in Sec. 4. Thus $W^{T_1}(\Sigma_{t_0})$ together with the map ${\hat \lambda^{t_0}}$ will be refer to as a local model of the {\bf first kind} for the local universal family of stable curves. 
  
  It is the precise version of the intuitive notion that the parameter $t\in W^{T_1}(\Sigma_0)$ is considered as a point near the ends of the stratum  $T_1$  representing 
  the stable curve $\Sigma_{t}=(S_t;  {\bf d}_t, {\bf x}_t )$ whose underlying  surface $S_t$ is deforming and degenerating along the ends, while its fixed part
  $K_t$  remains fixed such that  the relative locations of  ${\bf x}_t$ in $K_t$ are the same as the ones of ${\bf x}$ in $  K_0=K$. Note that in the  case used in \cite {lt}, the initial curve $\Sigma_0$ is "minimally'' stabilized. In this case, ${\bf x}_t$ is indeed {\bf  fixed } in the  above model.

   Now we define the required diffeomorphisms.
  In the case that $T_1=T_0$, the lowest stratum with $t_0=(b_0, 0)\in  W^{T_0}(\Sigma_{t_0})$,  $\lambda^{t_0}$ is just $\lambda^{b_0}=:\{\lambda^{b_0}_{b, v},  v\in T_0, b\in W^{T_0}(\Sigma_0)\}$
  defined   as follows.



  It is more convenient to define the inverse map  of $\lambda^{b_0}_{b, v}$.  For $v\in T_0, b\in W^{T_0}(\Sigma_0)$, $(\lambda^{b_0}_{b, v})^{-1}:(S_{b_0,v},  {\bf d}_{b_0, v}  )\rightarrow  (S_{b,v},  {\bf d}_{b, v}  ) $ is defined by 
  the following  conditions: 
  
  (i)  It is the "identity" map on the complement of the union of all  disks of radius $\epsilon_1$ centered
  at the double points  of the component $S_{b_0,v}$, denoted by $K_{b_0,v,\epsilon_1}$ under the identifications of
  ${\hat S}_{b}\simeq {\hat S}_{b_0}\simeq {\hat S}_{b=0}$  for $|b|$ and   $|b_0|$ sufficiently small. In fact under above identifications, we get the corresponding 
  identification $K_{b,v,\epsilon_1}\simeq K_{b_0,v,\epsilon_1}$ of the {\bf fixed } parts  given by  $\lambda^{b_0}_{b, v}$.

  (ii) On the disks $D_{\epsilon}(d_{b_0;vu})$ centered at the double point $d_{b;vu}$,
  It is the translation that brings $d_{b_0;vu}$ to $d_{b;vu}$, and $D_{\epsilon}(d_{b_0;vu})$ to $D_{\epsilon}(d_{b;vu})$. 
  
  (iii)  Under the above identification ${\hat S}_{b}\simeq {\hat S}_{b_0}$, the image $D_{\epsilon}(d_{b;vu})$ of 
  $D_{\epsilon}(d_{b_0;vu})$  can be considered   as a subset of  $D_{\epsilon_0}(d_{b_0;vu})\subset D_{\epsilon_1}(d_{b_0;vu})$ in $S_{b_0, v}$. Here of course, we assume that $\epsilon<\epsilon_0<\epsilon_1$ and $|b|$  and  $|b_0|$ sufficiently small.  Hence by (i) and (ii) above, using
  the above identification again, the  definition  for the rest of $(\lambda^{b_0}_{b, v})^{-1}$   is reduced to find a self  diffeomorphism of  $D_{\epsilon_1}(d_{b_0;vu})$ that is the identity map near the boundary extending
  the map already defined on $D_{\epsilon}(d_{b_0;vu})$ in (ii). This can be done
  by extending the corresponding vector field. Then the desired diffeomorphism is the time-$1$ map of the flow of the extended vector field.

  Note that the restriction map  $$\lambda^{b_0}_{b, v}: \cup_{[vu]\in E(T_0)}D_{\epsilon}(d_{b;vu})\cup  K_{b,v,\epsilon_1}\rightarrow \cup_{[vu]\in E(T_0)}D_{\epsilon}(d_{b_0;vu})\cup  K_{b_0,v,\epsilon_1}$$  is holomorphic and  preserves the "natural" metrics (Spheric, cylindrical or 'flat' ones).

 To move to the higher stratum $T_1$, recall that  for given $t=(b,a)\in W^{T_1}(\Sigma_0)$ with $t$ close $t_0=(b_0, a_0)$, at a
  double point $d_{b, vu}=d_{b, uv}$ where the gluing parameter $a_{vu}=a_{uv}\not = 0$,  the gluing of the pair of disks of radius $\epsilon $ on the components  $S_{b, v}$ and $S_{b, u}$ at the double points, $D_{\epsilon}(d_{b, vu})\#_{a_{vu}=a_{uv}}D_{\epsilon}(d_{b, uv})$ was defined in this section.
  
  The identification of the pair of disks $D_{\epsilon}(d_{b, vu})$ and $ D_{\epsilon}(d_{b, uv})$ with 
  $D_{\epsilon}(d_{b_0, vu})$ and $ D_{\epsilon}(d_{b_0, uv})$ 
  by $\lambda^{b_0}_{b, v}$ and 	$\lambda^{b_0}_{b, u}$ induces
  the corresponding  identification $$D_{\epsilon}(d_{b, vu})\#_{a_{vu}=a_{uv}}D_{\epsilon}(d_{b, uv})\simeq D_{\epsilon}(d_{b_0, vu})\#_{a_{vu}=a_{uv}}D_{\epsilon}(d_{b_0, uv}).$$  Applying this to each double point with $a_{uv}\not = 0, $
  We get  a family of identifications $S_{b_0,a}\simeq S_{b, a}$ smooth in $b$.
  
  Thus  the construction of $ \lambda^{t_0}_{t}$ with $t_0=(b_0, a_0)$, $t=(b,a)$ and $|t_0|$ and $|t|$ small can be  obtained by using  the following two families  of identifications of finite cylinders.
  
  (A)  When $l$ is close to $l_0$, there is a family of identifications
  $[-l_0, l_0]\times S^1\rightarrow [-l, l]\times S^1$ smooth in $l$ and induced from the corresponding identifications $[-l_0, l_0]\rightarrow [-l, l]$.
  We may assume that the restriction identifications to  $[-l_0=1, l_0-1]\times S^1$
  is the identity map.
  
  Applying these identifications, we get   a family of identifications
  $S_{b_0,a'}\simeq S_{b, a}$ with $|a'|=|a_0|$ and $arg \, a'=arg \, a$.

  (B) When $\theta\in S^1$ close to $1$, there is a smooth  $\theta$-dependent  family of identifications  $[-l_0, l_0]\times S^1\rightarrow [-l_0, l_0]\times S^1$ that is the  identity map  on $[-l_0, l_0-1]\times S^1$  and  is the rotation of angle $\theta$ on $\{l_0\}\times S^1.$ The effect of these identifications is to untwist the angular twisting in the gluing construction.
  Applying this to the identifications obtained above so far, we finally get
  the desire family of identifications $\lambda^{t_0}=\{\lambda^{t_0}_{ t}:S_t\rightarrow S_{t_0}\}$.
  
  Let $\epsilon_2>\epsilon_1$ and  assume that ${ D}_{\epsilon _1}({\bf d}_{b}\subset { D}_{\epsilon _2}({\bf d}_{b_0})$. Then define the ( "smallest")
  fixed part $K=:K_{\epsilon_2}$  to be the complement of  the union of all disks
  centered at double points of radius $\epsilon_2$ on $S_0$. Then we have the fixed parts
  $K_{t}\simeq K_{t_0}\simeq K_b\simeq K_{b_0}\simeq K$ in the corresponding 
  surfaces, on which the maps $\lambda^{t_0}_{ t}$ above are the "identity" map. 
  Note that the marked points ${\bf x}_t$ are lying on $K_{t}$.
  
  It follows from the construction above, the family of identification has the  properties  in the  above lemma.
  

 
 \subsection{Local universal family of the  second  kind}
 An important property the restriction  above local universal family ${\cal S}(\Sigma_0)$, a fixed stratum of type $T_1$ with 
 $T_1\geq T\geq T_0$, can be considered  as  a set of "moving" marked points ${\bf x}_t$ on  a fixed reference surface ${\hat S}_{t_0}$.   
  Then  the parameter $t\in \in W^{T_1}(\Sigma_{t_0})$, as part of the local coordinates of the moduli space ${\cal M}^{T_1}_{0,k}$ also describes the
  local moduli of the ("moving") distinguished points ${\bf p}_t$ on the ({\bf fixed }) reference curve $  {\hat S}_{t_0}$ of type $T_1$.

 
 The local universal family considered this way will be regarded as the second model.
 
 Thus in the second model, we are in the exact  the same situation as we were at the beginning of this section:  to parametrize stable curves ( considered as moving
 distinguished points on a fixed reference surface) in  a  fixed  stratum $T$, but with  the initial curve  $\Sigma_{t_0}$  instead of $\Sigma_0$.
 In particular in this model the   metric on  $  {\hat S}_{t_0}$ is fixed by the marking $  {\hat S}_{t_0}\simeq S^2$.
 

 
  To distinguish the situation here with the one before, instead of using the corresponding coordinate charts of the form $W^{T}(\Sigma_{t_0})$ with the local  parameter $b=b(t)\in W^{T}(\Sigma_{t_0})$  representing the moduli point of $\Sigma_t$, we will  denote $b=b(t)$ by $u:=u(t)$, and accordingly  ${\bf p}_{b(t)}$ by  ${\bf p}_{u}=:{\bf p}_{u(t)}$,
 $W^{T}(\Sigma_{t_0})$ by $ U^{T}(\Sigma_{u_0})$ with $u_0=u(t_0)$, etc.

 We state this more formally as following: the two  models of local universal family over the open set $N^{T}(\Sigma_{t_0})$ of ${\cal M}^{T}_{0, k}$ centered at $\Sigma_{t_0}=\Sigma_{u_0}$, denoted by 
  ${\cal S}|_{W^{T}(\Sigma(t_0))}$ and  $  {\cal S}|_{U^{T}(\Sigma(u_0))}$
   are identified by a fiber-wise analytic map $\phi^{-1}=\{\phi^{-1}_t\}$ with $\phi^{-1}_{t_0}:(S_{t_0},{\bf p}_{t_0})\rightarrow  (S_{u_0}; {\bf p}_{u_0})$  being the identification of the two (but the "same") central fibers.

 \begin{lemma}
 The  identification $\phi^{-1}_{t_0}:(S_{t_0},{\bf p}_{t_0})\rightarrow  (S_{u_0}; {\bf p}_{u_0})$ induces a family of component-wise biholomorphic identifications, denoted by 	 $\phi^{-1}_{t}: (S_{t},{\bf p}_{t}))\rightarrow  (S_{u}; {\bf p}_{u})$ with ${\hat S}_{u}$ being the fixed ${\hat S}_{u_0}$ such that the induced map $ {\underline \phi }^{-1}:W^{T}(\Sigma_{t_0})\rightarrow U^{T}(\Sigma_{u_0})$ 
 on the  coordinate charts given by  $t\rightarrow u$
   is a holomorphic identification. 
 Moreover  these maps $\{\phi^{-1}_t\}$  fit together to form an analytic   identification  between the corresponding universal families, denoted by $\phi^{-1}:{\cal S}|_{W^{T}(\Sigma_{t_0})}\rightarrow {\cal S}|_{U^{T}(\Sigma_{u_0})}$ that are holomorphic on the desingularizations  ${\hat{\cal  S}}|_{W^{T}(\Sigma_{t_0})}$ and $ {\hat {\cal S}}|_{U^{T}(\Sigma _{u_0})}$.

\end{lemma}

The above lemma is essentially a tautology in the situation above.  
However, it also follows from the local universal property of the family ${\cal S}$.






 
\section {Two types of neighborhoods near ends}

\subsection {"Base" deformations of a stable map $f=f_0:\Sigma_0\rightarrow M$}

\noindent
${\bullet }$ "Base" deformations of a initial  stable map $f=f_0:\Sigma_0\rightarrow M$ within the same stratum of type $T_0$.

\medskip

There are two kinds of  such "base" deformations of $f_0: \Sigma_0\rightarrow M$, denoted by  $\{{\tilde f}_b:\Sigma_b\rightarrow M, b \in W^{T_0}(\Sigma_0)\}$ and $\{f_b:\Sigma_b\rightarrow M, b \in W^{T_0}(\Sigma_0) \}$ respectively.

The first  one   is simply defined to be ${ \tilde f}_{b}=:\cup_{v\in T_0}{ \tilde f}_{b;v}$ with ${ \tilde f}_{b;v}=f_v$ under the assumption that $f$ is constant on
all the small disks near double points.

The last identity  makes sense as the domains of the two maps are the same $S^0_v$ after forgetting the distinguished points ${\bf p}_{b;v}$ and ${\bf p}_v^0$. 

To define the second deformation,
let  $\lambda_{b;v}^{b_0}:S_{b; v}\rightarrow S_v=:S_v^{b_0}$ in last section.
Then we define ${  f}_{b; v}:S_{b;v}\rightarrow M$ to be ${ f}_{b; v}=f_v\circ \lambda_{b;v}^{b_0}$
and ${ f}_b=\cup_{v\in T_0}{ f}_{b; v}:S_b\rightarrow M$.

Next we extend $\{f_b, b\in W^{T_0}(\Sigma_0)\}$ to a higher stratum $T_1>T_0$.

${\bf \bullet}$  "Base" deformation $\{f_t, t\in {\bar N}(\Sigma_0)\}.$

Given $f=f_{-}\bigvee f_{+};(D_{-}, d_{-})\bigvee_{d_{-}=d_{+}}(D_{+}, d_{+})\rightarrow M$ and a gluing parameter $a_0=exp\{-(s_0+t_0i)\}\not = 0$, to define  the gluing $\#_{a_0}(f_{-}, f_{+}): \#_{a_0}(D_{-}, D_+)\rightarrow M$ below, we introduce  the cylindrical coordinate $(s_{\pm}, t_{\pm})\in {\bf R}^{\pm}\times S^1	$ on $D_{\pm}$ by the 
identification of $D_{\pm}\simeq {\bf R}^{\pm}\times S^1$, $D_{\pm}$. Then $\#_{a_0}(D_{-}, D_+)$  is obtained by cutting of the part of $D_{\pm}$ with $|s_{\pm}|>-log|a_0|$ and glue the rest
 along the boundaries twisted with an angel $arg a_0.$ Thus $\#_{a_0}(D_{-}, D_+)\simeq [-log|a_0|, log|a_0|]\times S^1$ with the induced    cylindrical coordinate $(s, t)$ with $s=0$ corresponding to the "middle circle" and the other two 
  cylindrical coordinates $(s_{\pm}, t_{\pm})$ with $s_{\pm}=0$ corresponding to the two boundary circles.  Note that  $s_{\pm}=s{\mp}log|a_0|.$
 Then $\#_{a_0}(f_{-}, f_{+})$ is defined to be $\#_{a_0}(f_{-}, f_{+})(s,t)=exp_{f_{\pm}(d_{\pm})}(\beta_{-}(s){\hat f}_{-}(s, t)+(\beta_{+}(s){\hat f}_{+}(s, t)).$ Here ${\hat f}_{\pm}$ is a vector field over $D_{\pm}$ such that 
$f_{\pm}=exp_{f_{\pm}(d_{\pm})}	{\hat f}_{\pm}$, and $\beta_{\pm}$ are cut-off functions supported on $[-1,1]$ with $\beta_{-}+\beta_{+}=1.$

The deformation $f_t=(f_b)_a$ is then defined by implanting above construction to each double points $d_{vw}$ of $f_b$ with $a_{vw\not= 0.}.$

\subsection{Neighborhoods of the first kind }
Let $f_t, t\in {\bar W}(\Sigma_0)$ be the "base"   deformation of $f_0$ constructed above over   a  small coordinate/uniformizer ${\bar W}(\Sigma_0)$   of ${\bar N}(\Sigma_0)$ centered at $\Sigma_0$. 

A neighborhood $W^{\nu({\bf a})}_{\epsilon}(f_0, {\bf H}_{f_0})$ of $[f_0]$ in the space of unparametrized stable maps, as a slice, is defined to be a family of Banach manifolds over  the  "base"   deformation $\{f_t\},$
$$W^{\nu({\bf a})}_{\epsilon}(f_0, {\bf H}_{f_0})=\cup_{t\in  {\bar W}(\Sigma_0)}W_{\epsilon}^{ t; \nu(a)}(f_t, {\bf H}_{f_0}).$$

Here for each fixed $t$, 
$$W_{\epsilon}^{\nu(a) ;t}(f_t, {\bf H}_{f_0})
$$ 
$$ =\{h_t:(S_t, {\bf x}_t)\rightarrow 
(M, {\bf H}_{f_0})\,|\, \|h_t-f_t\|_{k, p; \nu(a)}<\epsilon\}.$$
Here $"a" $ is  gluing parameter in $t=(b,a)$.

\medskip
\noindent
${\bf \bullet }$ The $ \nu(a)$-exponentially weighted norm

 The norm $\|-\|_{k, p; \nu(a)}$  used in the definition above is  the $\nu$-exponentially weighted norm 
along the  all  "neck" areas  $N(b, a_{v_i, v_j})$ of $S_t$ with $a_{v_i, v_j}\not = 0 $ obtained by gluing from $S_0$ at    the double points $d_{v_i, v_j}(b)$; on the rest of $S_t$, the norm is just
the usual $L_k^p$-norm.  More specifically, recall that "neck" areas  of   $S_t$ with $a_{v_i, v_j}\not = 0 $ and $t=(b, a)$,  $N(b, a_{v_i, v_j})\simeq (-|log |a_{v_i, v_j}||, |log |a_{v_i, v_j}||)\times S^1$.  The  weight
 function  $\nu(a_{i, j})$  is  equal to $exp\{\nu |s_{\pm}|\}$ for points in $N(b, a_{v_i, v_j})$ with 
 $s_{\pm}\in (0, |log |a_{v_i, v_j}|\,|-2)$ and  is a smooth function equal to  1 on two ends of the neck. Outside these necks, $\nu(a)=1$ so that these weight $\nu(a_{i, j})$
 functions   together defines a smooth   weight function $\nu(a)$.
 Here $\nu$ is  a fixed  positive constant with $\nu< (p-2)/p.$
 The $L^p_{k, \nu(a)}$-norm  $\|h_t\|_{k, p, \nu(a)}$ then is  defined to be 
 $\|\nu(a)\cdot h_t\|_{k, p}.$ 
 
 For each fixed $t$, the norm so defined makes $W_{\epsilon}^{\nu(a) ;t}(f_t, {\bf H}_{f_0})
 $   became a  Banach manifold.  However, on  $W^{\nu({\bf a})}_{\epsilon}(f_0, {\bf H}_{f_0})$, this $t$ -dependent family of  norms  does  not necessarily define 
 a topology without further conditions
 as in [L6?]. The reason for this is that  the  mixed $L_{k}^p$/$L_{k, \nu(a)}^p$-norm   used  here is  not continuous when $h_t$ is moving from higher stratum to the lower ones.
 On the other hand, on each fixed stratum, near any given point the $L_{k}^p$ and $L_{k, \nu(a)}^p$ norms are "locally" equivalent so that the resulting space is at least  a (topological ) Banach manifold.

 More specifically,  consider the decomposition of $W^{\nu({\bf a})}_{\epsilon}(f_0, {\bf H}_{f_0})$  into its open strata,
$$W^{\nu({\bf a})}_{\epsilon}(f_0, {\bf H}_{f_0})=: \cup_{T\geq T_0}W^{ \nu (a_T),T}_{\epsilon}(f_0, {\bf H}_{f_0})$$ with each  $$W^{ \nu (a_T), T}_{\epsilon}(f_0, {\bf H}_{f_0})=:\cup_{t\in W^{T}_{\epsilon}(\Sigma_0)}W_{\epsilon}^{ \nu(a),t}(f_t, {\bf H}_{f_0}).$$ Then the collection of all such neighborhoods $W^{ \nu (a_T), T}_{\epsilon}(f, {\bf H}_{f_0})$ of a fixed stratum $T$ generate a topology by the above mentioned local equivalent equivalence of the $L_{k, \nu(a)}^p$-norm with standard  $L_k^p$-norm.

Since in order to prove the main results of this paper on the higher smoothness of the admissible perturbations, we need to localize further by  using  small  neighborhoods of 
 any   given element $g_{t_0}$ in $W^{ \nu (a_T), T}_{\epsilon}(f, {\bf H}_{f_0})$. We need spell out more on the existence of such neighborhoods.





${\bullet}$  Neighborhoods of the   first  kind  on  $W^{\nu({\bf a})}_{\epsilon}(f_0, {\bf H}_{f_0})$:

Fix a stratum of type $T=T_1\geq T_0.$ Recall that $T_0$ is the tree associated with the lowest stratum that $f_0$ lies on.
Give $\{g_{t_0}:(S_{t_0}, {\bf x}_{t_0})\rightarrow (M, {\bf H}_{f_0})\}\in 
 W_{\epsilon}^{\nu (a_{T_{1}}), T_{1} }(f_0, {\bf H}_{f_0})$ with $t_0=(b_0, a_0)\in W^{T_1}(\Sigma_0)$ and  $\Sigma_{t_0}=(S_{t_0}, {\bf p}_{t_0})$,
	 a neighborhood of $g_{t_0}$ of the {\bf first kind} in 
  $ W_{\epsilon}^{ \nu (a_{T_{1}}),T_{1}}(f_0, {\bf H}_{f_0})$,  denoted by

 $$W^{\nu(a_{T_{1}}), T_{1}}_{\epsilon '}(g_{t_0}, {\bf H}_{f_0})=:W^{\nu(a_{T_{1}}), T_{1}}_{\epsilon '}(g_{t_0}, {\bf H}_{f_0}; g_t)=\{h_t:S_t\rightarrow M|\,\, \|h_t-g_t\|_{k, p, {\nu(a)}}<\epsilon'\}$$   with  ${\epsilon '}<<\epsilon ,$
 can be defined in  a few equivalent ways that we describe now.
  
   Here $g_t:S_t\rightarrow M$  is the "base" deformation of $g_{t_0}$ inside  $  W_{\epsilon}^{ \nu (a_{T_{1}}),T_{1}}(f_0, {\bf H}_{f_0})$, 
 similar   to the initial  deformation $f_t$.    Since   topological type  of $S_t$  is  fixed, 
   we require that 
   the deformation has a form 
   $g_t=g_{t_0}\circ T^{t_0}_t$,  where $$T=:T^{t_0}=\{T_t^{t_0}\}: {\cal S }|_{W^{T_1}_{\epsilon'}(\Sigma_{t_0})}=\cup_{t\in {W^{T_1}_{\epsilon'}}(\Sigma_{t_0})} S_t \rightarrow  S_{t_0}$$
   is  a smooth  family of diffeomorphisms.  
   
    Thus we need   establish the existence of the required deformation $g_t$. 
     The key step is the following lemma.

\begin{lemma}
	Given the base deformation $\{f_t\}$ of $f_0$ defined earlier in this section,  fix a member 
	$f_{t_0}:S_{t_0}\rightarrow M$  of type $T_1$ in the deformation, there is a  smooth family of diffeomorphisms $T=:T^{t_0}=\{T_t^{t_0}\}: {\cal S }|_{W^{T_1}_{\epsilon'}(\Sigma_{t_0})}\rightarrow S_{t_0}$ such that $\lim_{t\rightarrow t_0}\|f_t-f_{t_0}\circ T^{t_0}_{t}\|_{k, p, \nu(a)}=0$ in an uniform manner in $t_0$ for $t_0$ varying in a compact set.
\end{lemma}
  
 Above lemma implies the following two lemmas 
 
 \begin{lemma}
 	Given $g_{t_0}\in W_{\epsilon}^{ \nu (a_{T_{1}}),T_{1}}(f_0, {\bf H}_{f_0})$. There exists a  "base" deformation  $g_t$  inside  $  W_{\epsilon}^{ \nu (a_{T_{1}}),T_{1}}(f_0, {\bf H}_{f_0})$ for $|t-t_0|<<\epsilon$ with   
 	$g_t=g_{t_0}\circ T^{t_0}_t$.

 \end{lemma}
 
 \begin{lemma}
 	These $W^{\nu(a_{T_{1}}), T_{1}}_{\epsilon '}(g_{t_0}, {\bf H}_{f_0})$  generate a topology on $ W_{\epsilon}^{ \nu (a_{T_{1}}),T_{1}}(f_0, {\bf H}_{f_0})$.
 \end{lemma}
 
 We note that the required $T^{t_0}_t$ can be taken as the particular family of diffeomorphisms $\lambda^{t_0}_t$ defined before.
 Since our main concerns of this paper is the stratified smoothness of generic perturbations, we will not give the proofs of above lemmas. They will be given in [L?].

 The neighborhoods  $W_{\epsilon'}^{\nu(a_{T_{1}}), T_{1}}(g_{t_0}, {\bf H}_{f_0})$   here will be call the ones of the {\bf first} kind.

 Here are some variations or  related constructions:
 
 (1) In the definition of $W^{\nu(a_{T_{1}}), T_{1}}_{\epsilon '}(g_{t_0}, {\bf H}_{f_0})$, the  conditions that $|t-t_0|< \epsilon'$ and $\|h_t-g_t\|_{k, p; \nu(a)}< \epsilon'$  can be replaced by $|t-t_0|+\|h_t-g_t\|_{k, p; \nu(a)}<\epsilon'.$

 (2)   In the above definition of   the $L_{k, \nu(a)}^p$-norm, a $t$-dependent metric $m_t$　 on the domain $S_t$.
 is used.  Using the diffeomorphisms $T_{t}^{t_0}:S_t\rightarrow S_{t_0}$ 
 to pull-back the fixed metric $m_{t_0}$, we  get a family of metrics $ (T_{t}^{t_0})^*(m_{t_0})$ on $S_t$ and the corresponding  $L_{k, \nu(a)}^p$-norms and  neighborhoods $W^{\nu(a_{T_{1}}), T_{1}}_{\epsilon '}(g_{t_0}, {\bf H}_{f_0})$. Since these families of the metrics are uniformly equivalent for
 $|t-t_0|\leq \epsilon''$, the resulting neighborhoods defined this way  are equivalent to the previous  one.
 Thus, upto the effect of $T_t^{t_0}$, we can use a fixed reference metric on $m_{t_0}$ to define the norm.
 
 This implies that the second  type of the neighborhoods  defined below  is (topologically) equivalent to the ones  above.

 (3) Consider  the deformations 
 of  $h_t=h_{t_0}\circ T_t^{t_0}$ for all $h_{t_0}:S_0\rightarrow M $ with  $h_{t_0}$ in the central slice $ W^{\nu(a_{0}),t_0,  T_{1}}_{\epsilon '}(g_{t_0}, {\bf H}_{f_0})$.  Denote the collection of such $h_t$ with $|t-t_0|<\epsilon'_1$ by 
 $W'^{\nu(a_{T_{1}}), T_{1}}_{\epsilon '_1}(g_{t_0}, {\bf H}_{f_0})$.
 
 \begin{lemma}
 	The neighborhoods $W'^{\nu(a_{T_{1}}), T_{1}}_{\epsilon '_1}(g_{t_0}, {\bf H}_{f_0})$
 	so defined are equivalent to the ones defined before.
 \end{lemma}

 For the proof of this lemma we refer to [L?] again.
  \medskip
  \subsection {   Neighborhoods
  of the second kind}

  Recall that   there is a family of 
  biholomorphic identifications	 $\phi^{-1}_{t}: (S_{t},{\bf p}_{t})\rightarrow  (S_{u}; {\bf p}_{u})$ which transforms the family of varying curves $\{(S_{t},{\bf p}_t)\}$ with {\bf fixed}  $(K_t,{\bf x}_t)=(K_0,{\bf x}_0)$ with
  the curve with {\bf fixed} components  but with a family of varying distinguished  points,   $\{({\hat S}_{u_0}, {\bf p}_u)\}$.

It induces a bijection  $$\Phi= \cup_{t}\Phi_t: W_{\epsilon'}^{\nu(a_{T_1}), T_{1}}(g_{t_0}, {\bf H}_{f_0})=:\cup_{t\in W^{T_1}(\Sigma_{t_0})} W_{\epsilon'}^{\nu(a), t}(g_{t_0}, {\bf H}_{f_0}) $$ $$  \rightarrow  U_{\epsilon'}^{ T_{{1}}}(g_{u_0}, {\bf H}_{f_0})=\cup_{u\in U^{T_1}(\Sigma_{u_0})}
  U^{u}_{\epsilon'}(g_{u_0}, {\bf H}_{f_0}) $$ by pull-backs given by $\Phi_t(h_t)=
  h_{t}\circ \phi_{t}$ denoted by $h_u.$ 
  
    The subspace  $U^{u}_{\epsilon'}(g_{u_0}, {\bf H}_{f_0}) $  here  consists those  stable maps $h_u:(S_u, {\bf x}_u)\rightarrow
  (M, {\bf H}_{f_0})$ with {\bf fixed} marked points ${\bf x}_u$ on  the  fixed  domain ${\hat S_u}={\hat S_{u_0}}$  such that $\|h_u-g_u\|_{k, p}<{\epsilon'}$. Here we give
  each component of  ${ \hat S}_{u_0}$ the spherical metric and use it to defined the  above  $L_k^p$-norm, and $\{g_u=g_t\circ \phi_t\}$ is the transformed base family by $\{\phi_t\} $.
  
  Since for $|t-t_0|<\epsilon'$, the norms of the above  two spaces are equivalent, $\Phi$ is a homeomorphism to its image.
  
   The neighborhoods $U_{\epsilon'}^{ T_{{1}}}(g_{u_0}, {\bf H}_{f_0})$  will be refereed as  of the {\bf second} kind. 
   By the remark (2) above, the neighborhoods  here are equivalent to the ones of the first kind above.

Note that the lowest stratum, $W^{\nu(a_{T_0}, T_0 )}_{\epsilon}(f_0, {\bf H}_{f_0})$
 is the same as  $U^{T_0}_{\epsilon}(f_{u=0}, {\bf H}_{f_0})$ since the norm used here is just the usual $L_k^p$-norm without exponential  weight.

It is easy to see that   
$\Phi_t$ is a diffeomorphism. 

\begin{lemma}
	For each fixed $t$ and hence $u$, $$\Phi_t: W_{\epsilon'}^{\nu(a), t}(g_{t_0}, {\bf H}_{f_0})  $$ $$  \rightarrow U^{u}_{\epsilon'}(g_{u_0}, {\bf H}_{f_0})
 $$ is a	diffeomorphism.
\end{lemma}

\proof

The only difference between $U_{\epsilon'}^{ T_1}(g_{u_0}, {\bf H}_{f_0})$ and $W_{\epsilon'}^{\nu(a_{T_1}), T_1}(g_{t_0}, {\bf H}_{f_0})$ is that the domains of stable maps in $W_{\epsilon'}^{\nu(a_{T_1}), T_1}(g_{t_0}, {\bf H}_{f_0})$ are varying depending on the parameter $t$  while the domains of  elements in $U^{ T_{1}}(g({t_0}), {\bf H}_{f_0})$ are the fixed  ${\hat S}_{u_0}$ but with moving  distinguished  points  parametrized by $u=u(t)$.  When $t$ is fixed, the elements
in $U_{\epsilon'}^{u}(g_{u_0}, {\bf H}_{f_0})$ or $W_{\epsilon'}^{\nu(a), t}(g_{t_0}, {\bf H}_{f_0})$ have the same  domain with fixed components  and distinguished  points under the  identification map $\phi_t$. Moreover the Banach norm on these two  spaces are equivalent. Hence  the induced map $\Phi_{t}$ is a diffeomorphism.

\QED

From now on, if there is no confusion, we will drop the subscript $\epsilon$ that describes the size of a neighborhood.

Next we define  the neighborhoods that are still  second type obtained from $U^{ T_{1}}(g_{u_0}, {\bf H}_{f_0})$ by 
dropping some marked points. 

Given $\{h:S_u\rightarrow M\} \in U^{ T_{1}}(g_{u_0}, {\bf H}_{f_0})$ with fixed ${\hat S_u}={\hat S}_{u_0}$, among the $k$ marked  points ${\bf x}_u$ on ${\hat S}_{u_0}$ we    select $m$ points, denoted by ${\bf x}^r_u$ such that $ (S_u,{\bf x}^r_u) $ is still stable. Denote the resulting stable curve  with $m$ marked points
by   $\Sigma^r_u=:(S^r_u, {\bf x}^r_u)$  with fixed ${\hat S}^r_u={\hat S}^r_{u_0}$, which is the same as ${\hat S}_{u_0}$ as a surface.

Then the  map  $ h^r:(S^r_u, {\bf x}^r_u)\rightarrow (M, {\bf H}^r_{f_0})$ is defined to be the same  map  $ h$ as before but forgetting  the rest of the marked points, denoted by  ${\bf x}^c_u$,  in ${\bf x}_u$, where ${\bf H}^r_{f_0}$ is the corresponding
selection of local hypersurfaces.  Note that $h^r_u({\bf x}^c_u)\in {\bf H}^c_{f_0}.$

Then the collection of all such $ h^r$  obtained from  $h\in U^{ T_{1}}(g_{u_0}, {\bf H}_{f_0})$  by dropping $k-m$ marked points will be denoted by   $U^{ T_{1}}(g^r_{u_0}, {\bf H}^r_{f_0})$ with the  centered 
$ g^r_{u_0}$.    

The process above of course depends on the following choices: (i) the selection of ${\bf x}^r_u$, (ii) an order for ${\bf p}_u\subset {\hat S}_u={\hat S}_{u_0}$  that induces an order for  ${\bf p}^r_u\subset S^r_u={\hat S^r}_{u_0}$.
In the following,  we fix one of such  choices, labeled by the superscript $r$ in the  notations here. For each $v\in T_1$, by identifying the first three points in $({\bf p}^r_u)_v$ and  $({\bf p}_u)_v$,  we get the holomorphic identifications $\psi_v:({\hat S^r}_u)_v=({\hat S^r}_{u_0})_v\rightarrow ({\hat S}_u)_v=({\hat S}_{u_0})_v  $ and
$\psi=\cup_{v \in T_1}\psi_{v}:{\hat S}^r_{u_0}\rightarrow {\hat S}_{u_0}, $  which induces the above drop-marking map, denoted by $$\Psi =\cup_{ v\in T_1}\Psi_v:U^{ T_{1}}(g_{u_0}, {\bf H}_{f_0})\rightarrow U^{ T_{1}}(g^r_{u_0}, {\bf H}^r_{f_0})$$  given by 
pull-back by $\psi$, 
$h\rightarrow h\circ \psi$ denoted by $h^r $.

In next section we  will show that this map  is a    diffeomorphism.

Note that the identifications  $\psi_v:({\hat S^r}_u)_v=({\hat S^r}_{u_0})_v\rightarrow ({\hat S}_u)_v=({\hat S}_{u_0})_v  $
 above also give  canonical identifications of these surfaces with the fixed $(S^2,; 0, 1, \infty)$  so that they become "marked"  surfaces.

\section{Stratified smooth structures on  neighborhoods}

 The smooth structures on $U^{ T_{1}}(g^r_{u_0}, {\bf H}^r_{f_0})$  and  $U^{ T_{1}}(g_{u_0}, {\bf H}_{f_0})$  are defined similarly,
 obtained  as  $C^{m_0}$ submanifolds $E_m^{-1}({\bf H}^r_{f_0})$ of ${\widetilde U}^{ T_{1}}(g_{u_0})\times {\hat S}_{u_0}^{m }$  and  $E_k^{-1}({\bf H}_{f_0})$ of  ${\widetilde U}^{ T_{1}}(g_{u_0})\times {\hat S}_{u_0}^{k }$ respectively. Here $E_l:{\widetilde U}^{ T_{1}}(g_{u_0})\times {\hat S}_{u_0}^{l }\rightarrow M^l$  is the $l$-fold  total evaluation map at $l$ selected
 marked points among the $k$ marked points on the domain $( {\hat S}_{u}, {\bf x}_u)=({\hat S}_{u_0}, {\bf x}_u)$. 
 Since  $E_l$ is of class $C^{m_0}$ (see  \cite {L8} ), it is easy to see that
  it is a $C^{m_0}$-submersion so that above two subsets have   $C^{m_0}$-smooth structures.

\begin{pro}
$\Psi :U^{ T_{1}}(g_{u_0}, {\bf H}_{f_0})\rightarrow U^{ T_{1}}(g^r_{u_0}, {\bf H}^r_{f_0})$  is a local diffeomorphism at $g_{u_0}$. 
\end{pro}

\proof


Note that $\Psi :U^{ T_{1}}(g_{u_0}, {\bf H}_{f_0})\rightarrow U^{ T_{1}}(g^r_{u_0}, {\bf H}^r_{f_0})$  is a bijection.
 Indeed since $g^r_{u_0}({\bf x}^c_{u_0})\in {\bf H}^c_{f_0}$ and for any $h\in  U^{ T_{1}}(g^r_{u_0}, {\bf H}^r_{f_0})$, we already have $h({\bf x}^r_u)\in {\bf H}^r_{f_0}$,  when  $U^{ T_{1}}(g^r_{u_0}, {\bf H}^r_{f_0})$ is small enough and ${\bf x}^r_u$ as above is fixed, by implicit function theorem, the equation on ${\bf x}^c_u$,  $h({\bf x}^r_u, {\bf x}^c_u)\in {\bf H}^c_{f_0}$ has an unique  solution that is close to ${\bf x}^c_{u_0}$. This  proves  that $\Psi=:\Psi^r$ is a bijection with $\Psi^{-1} $ sending $\{h^r:({\hat S}^r_u={\hat S^r}_{u_0},  {\bf x}^r_u)\rightarrow M \}$ to  $\{h:({\hat S}_u={\hat S}_{u_0},  {\bf x}^r_u,  {\bf x}^c_u )\rightarrow M \}$. 
 
   Thus the map $\Psi$ is the restriction to a $C^{m_0}$ submanifold of the obvious smooth  projection  $\pi:{\widetilde U}^{ T_{1}}(g_{u_0})\times {\hat S}_{u_0}^{k}\rightarrow {\widetilde U}^{ T_{1}}(g_{u_0})\times {\hat S}_{u_0}^{m }$ sending 
  $k$ marked points to the corresponding $m$ marked points.

 \QED 
 
 {\bf  Note}:  This process of dropping-adding marked points given by  $\Psi$  is used in \cite {5} by  requiring  that domains of the elements in $W^{ T_1}(g^r_{u_0}, {\bf H}^r_{f_0})$ are  stabilized with minimal number of marked points.
 

 Next we  defined a  smooth structure centered at $g_{t_0}$
for  $W^{\nu(a_{T_1}), T_1}(g_{t_0}, {\bf H}_{f_0}) $
  or $W^{\nu(a_{T_1}), T_1}(g_{t_0}) $. This can be done by using the  family of smooth
  identifications $\lambda^{t_0}=\{\lambda^{t_0}_{ t}:\Sigma_t\rightarrow \Sigma_{t_0},  t\in W^{T_1}(\Sigma_{t_0})\}$ that is the identity map on the "small fixed part" $K_{t,\epsilon_2}({\bf d}(b_0))$ "centered" at ${\bf d}(b_0)$ defined  before in Sec. 2.  These maps give  rise a  smooth trivialization   of  the local universal family $({\cal S}|_{W^{T_1}(\Sigma_{t_0})}\rightarrow  W^{T_1}(\Sigma_{t_0}))\simeq  \Sigma_{t_0}\times W^{T_1}(\Sigma_{t_0})$.

  Now the smooth structure on   $W^{\nu(a_{T_1}), T_1}(g_{t_0}, {\bf H}_{f_0}) $ can be defined by the identification $\Lambda^{t_0}:  W^{\nu(a_0), t_0}(g_{t_0}, {\bf H}_{f_0}) \times N^{T_1}(\Sigma_{t_0}) \rightarrow W^{\nu(a), T_1}(g_{t_0}, {\bf H}_{f_0})$ defined by $\Lambda^{t_0}(h, t)=h\circ \lambda^{t_0}_t.$ Note that the norms are equivalent under the map $\Lambda^{t_0}$ by the remark/note  (3) in last section.
  This gives a  smooth structure on each $W^{\nu(a), T_1}(g_{t_0}, {\bf H}_{f_0}) .$ Of course transition functions between two such neighborhoods of type $T_1$
  are  only  continuous.  The end  $W^{\nu(a_{T_1}), T_1}(f_{0}, {\bf H}_{f_0}) $ then is covered by such  neighborhoods.

  \section{Higher-degree stratified smoothness of \\ the perturbations }


  In this section we  give a proof that  the  ${\xi}$ on 
  $W^{\nu(a_0), T_1}(g_{t_0}, {\bf H}_{f_0})$ defined below is of class $C^{m_0}$ viewed in any other local slices.


   We assume that ${\xi}$ is obtained from ${\xi}^K=\oplus_{v\in T_0}{\xi}^{K_v}$  defined below.

   Here 
     $K=K_0=\cup_{v\in T_0}K_v\subset S_0$ is the  fixed part lying on the  initial curve $S_0$ with subscript "$0$" corresponding to $t=(a, b)=(0, 0)$.  

  For each $v\in T_0$,  let  ${\widetilde W}(f_{K_v})$ be  the collection of $L_k^p$-maps $g=g_v:K_v\rightarrow M$ such that $\|g-f_{K_v}||_{k, p}<\epsilon$. Here $f_{K_v}=f|_{K_v} $,  the restriction of the initial map $f$ to $K_v$. We give
  $K_v$ the induced metric from $S_v$. The  bundle $({\cal L}^{K_v}, {\widetilde W}(f_{K_v}))$ is defined as following: for any $h_v\in {\widetilde W}(f_{K_v})$, the fiber ${\cal L}^{K_v}|_{h_v}=(L_{k-1}^p)_0(K_v,h_v^*(E) )$ consists of all $L_{k-1}^p$-sections with compact support in the interior of $K_v$.   
 
  Let ${\widetilde W}(f_{K})=\prod_{v\in T_0}{\widetilde W}(f_{K_v})$ and ${\cal L}^{K}
  =\oplus_{v\in T_0}{\cal L}^{K_v}$. 
  
 Now for each $v\in T_0$, fix a section ${\xi}^{K_v}:{\widetilde W}(f_{K_v})\rightarrow {\cal L}^{K_v}$ of class $C^{\infty}$ satisfying the conditions $C_1$ and $C_2$. Then the section  ${\xi}^{K}=:\oplus_{v\in T_0}{\xi}^{K_v}:{\widetilde W}(f_{K})\rightarrow {\cal L}^{K}.$

  By  the identification of $K_t\simeq K=K_0$, for any fixed
   $t\in W^{T_1}(\Sigma_{t_0})$, and $g_t\in W^{\nu(a), t}(g_{t_0}, {\bf H}_{f_0})\subset W^{\nu(a_{T_1}), T_1}(g_{t_0}, {\bf H}_{f_0}) $, 
   we get the induced section $\xi^t: W^{\nu(a), t}(g_{t_0}, {\bf H}_{f_0})\rightarrow {\cal L}^t$ defined by $\xi^t(g_t)=:\xi^K(g_t)|_{K_t}.$
   It is easy to see  that the standard local trivializations for  
   the bundle $({\cal L}^{K}\rightarrow {\widetilde W}(f_{K}))$ and $({\cal L}^{t}\rightarrow {\widetilde W}^{\nu(a), t}(g_{t_0}, {\bf H}_{f_0}))$ are  compatible with respect to the above identifications
    of $K_t\simeq K_0$ so that for each fixed $t$, $\xi^t$ is still smooth
    and satisfies the condition $C_1$ and $C_2$.

    It follows 
    that 
 ${\xi}^{K}$ becomes
    a section ${\xi}$ on $W^{\nu(a_{T_1}), T_1}(g_{t_0}, {\bf H}_{f_0}) $,  defined  by ${\xi}=\cup_{t\in W(\Sigma_{t_0})}\xi^t$.

     In fact  it becomes a section
        ${\widetilde\xi}$ on the larger space ${\widetilde W}^{\nu(a_0), T_1}(g_{t_0} )$ without the constraints given by ${\bf H}_{f_0} $.  

         
     Recall that these  identifications $\lambda^{t_0}=:\{\lambda^{t_0}_t:S_t\rightarrow S_{t_0}\}$, 
     ${\cal S}|_{W^{T_1}(\Sigma_{t_0})}\rightarrow S_{t_0}$,
     	 induce a product structure $${\widetilde  W}^{\nu(a_{T_1}), T_1}(g_{t_0} )\simeq {\widetilde W}^{\nu(a_0), t_0}(g_{t_0} ) \times W^{T_1}(\Sigma_{t_0}) ,$$  which in turn gives a smooth structure on  ${\widetilde  W}^{\nu(a_{T_1}), T_1}(g_{t_0} )$.
     Moreover, for $t$ is sufficiently close to $t_0$ the map $\lambda_{t}^{t_0}:S_t\rightarrow S_{t_0}$ induces an identification of 
     ${\cal L}^t_{g_t}$ with ${\cal L}^{t_0}_{g_{t_0}}$ for 
     ${\hat g}_t={\hat g}_{t_0}\circ (\lambda_{t}^{t{_0}})$ essentially by the pull-back of $\lambda_{t}^{t_0}$. Indeed, in the case, that $E=TM$ and  ${\cal L}^t_{g_t}=L_{k-1}^p(S_t, {g_t}^*(E))$ or  ${\cal L}^t_{g_t}=L_{k-1}^p(S_t, {g_t}^*(E)\otimes \Lambda^ {1})$, it is exactly
      given by the pull-back. 
     
   For   ${\cal L}^t_{g_t}=L_{k-1}^p(S_t, {g_t}^*(E)\otimes \Lambda^{0, 1})$  and  ${\cal L}^{t_0}_{g_{t_0}}=L_{k-1}^p(S_t, {g_{t_0}}^*(E)\otimes \Lambda^{0, 1}),$
     the identification is given by composition  of the pull-back by $\lambda_{t}^{t_0}$ with the map induced by the projection $  \Lambda^{ 1}\rightarrow \Lambda^{0, 1}$ since  $\lambda_{t}^{t_0}$ is not holomorphic away from the fixed part $K$ as already observed in \cite {lt}.
     Combing this  with the 
     local trivializations of the bundles $({\cal L}^t\rightarrow {\widetilde W}^{\nu(a_0), t}(g_{t}))$ above, this gives a
     trivialization of 
     the bundle $\{{\cal L}^{T_1}=\cup_{t\in W^{T_1}(\Sigma_{t_0})}{\cal L}^t
     	\rightarrow {\widetilde W}^{\nu(a_{T_1}), T_1}(g_{t_0})=
     	\cup_{t\in W^{T_1}(\Sigma_{t_0} )}{ \widetilde W}^{\nu(a_0), t}(g_{t})\}$  centered at $g_{t_0}$. 
     	
     	Using the  fact that $\lambda_{t}^{t_0}$ is just the identity map on $K=K_{t}=K_{t_0}$ (and hence holomorphic) it is easy to check that with respect to this product smooth structure and local trivialization,  $\xi /{\tilde  \xi}$ so defined   is 
     	"constant" along $W^{T_1}(\Sigma_{t_0}) $-directions in the sense that
     	for $(h, t)\in { W}^{\nu(a_0), t_0}(g_{t_0} ) \times  W^{T_1}(\Sigma_{t_0})\simeq { W}^{\nu(a_{T_1}), T_1}(g_{t_0} )  ,$
     	$\xi(h, t)=\xi^{t_0}(h)(=\xi^K(h|_K))$ ( similarly for ${\hat \xi}$). Hence it is   a smooth section.
      Clearly the condition $C_1$ and $C_2$ still hold for ${\tilde \xi}$. 
      
      The discussion here can be reformulated in the lemma.
      
    \begin{lemma} 
    On ${ \widetilde W}^{\nu(a_{T_1}), T_1}(g_{t_0})$,
     ${\tilde  \xi}=(\lambda^{t_0})^*({\tilde  \xi}^{t_0}).$ 
      Here the section ${\tilde  \xi}|_{ \widetilde W^{\nu(a_{0}), t_0}(g_{t_0})}$ along the central slice ${ \widetilde W^{\nu(a_{0}), t_0}(g_{t_0})}$ is denoted by  ${\tilde  \xi}^{t_0}$.  
    
       \end{lemma} 
      
    \proof  
      
 This essentially is a tautology.   
   \QED   
   
   Thus ${\tilde \xi}$ is the $G_e^1$-extension of ${\tilde {\xi}}^{t_0}$. Here $G_e^1=W^{T_1}(\Sigma_{t_0})$ considers as a family of diffeomorphisms $\{\lambda^{t_0}_t:S_t\rightarrow S_{t_0}\}$ parametrized by $t\in W^{T_1}(\Sigma_{t_0})$.   This interpretation proves the above lemma again. Similar interpretations using  $G_e^i$-extensions with $i>1$  prove the two main theorems
    below.

     The first main theorem  is the following.
    
    
    \begin{theorem}
    	The section ${\tilde  \xi}^{\Phi}$ is of class $C^{m_0}$.
    	  Consequently its
    	 restriction ${ \xi}^{\Phi}$  to the $C^{m_0}$-submanifold ${ U}^{  T_{{1}}}(g_{u_0}, {\bf H}_{f_0}) $ with respect to the second smooth structure is of class $C^{m_0}$. Here  ${\xi}^{\Phi}$ is obtained
    	 from the section  ${\xi}$ on 
    	 $W^{\nu(a_0), T_1}(g_{t_0}, {\bf H}_{f_0})$ by the transformation $\Phi$, similarly for ${\tilde  \xi}^{\Phi}$. 
    	 The section  ${ \xi}^{\Phi, \Psi}$ on 
    	  $U^{ T_{1}}(g^r_{u_0}, {\bf H}^r_{f_0})$ 
    	   induced by the diffeomorphism $\Psi :U^{ T_{1}}(g_{u_0}, {\bf H}_{f_0})\rightarrow U^{ T_{1}}(g^r_{u_0}, {\bf H}^r_{f_0})$  is of class $C^{m_0}$ as well.  
 \end{theorem} 
 
 \proof
 
 The first statement  follows from the fact that ${\tilde \xi}^{\Phi}$ is the $G_e^2$-extension of ${\xi}^{t_0}$ by pull-backs of the elements in $G_e^2$.  Here $G_e^2=W^{T_1}(\Sigma_{t_0})(=U^{T_1}(\Sigma_{u(t_0)}))$ considered as the 
 family of diffeomorphisms $\{\lambda^{t_0}_t\circ \phi_t:S_{u(t)}\rightarrow S_t\rightarrow S_{t_0}\}$ parametrized by $t\in W^{T_1}(\Sigma_{t_0})$. 
 Since ${\xi}^{t_0}$ satisfies the conditions $C_1$ and $C_2$, a obvious generalization of the main theorem in the first paper of these sequel implies that the $C^{m_0}$-smoothness of the $G_e^2$-extension of ${\xi}^{t_0}$ above so that 
  ${\tilde \xi}^{\Phi}$ is of class $C^{m_0}$.

 To prove the last statement, we note that
  in addition to diffeomorphism $\Psi:U^{ T_{1}}(g_{u_0}, {\bf H}_{f_0})\rightarrow U^{ T_{1}}(g^r_{u_0}, {\bf H}^r_{f_0})$ of that identified the  bases, the standard local trivializations  of the bundles ${\cal L}$ centered at $ g_{u_0}$ and $g^r_{u_0}$  using parallel transport along shortest connecting  geodesics are also the "same"
   in the sense that  the trivialization for ${\cal L}$ on  $U^{ T_{1}}(g^r_{u_0}, {\bf H}^r_{f_0})$ automatically give  the one for  ${\cal L}$ on the other by  our assumption that each local hypersurfaces in ${\bf H}_{f_0}$ is geodesic submanifold.
  
 \QED

 The idea above can be used  to proof the smoothness of ${\xi}$ viewed in any  other chart $W'^{T_1}(f_0', {\bf H}')$:  that is  to define the corresponding $G_e^i$-extension of the same section ${\xi}^{t_0}$.
 
 We start with the neighborhood $U'^{T_1}(g'_{u'_0}, {\bf H}')$ with the class $[g'_{u'_0}]=[g_{u_0}]$. Here $g_{u_0}=g_{u(t_0)}:(S_{u_0}, {\bf x}(u))\rightarrow (M,{\bf H}_{f_0})$ is the center of $U^{T_1}(g_{u_0}, {\bf H}_{f_0})$, similarly for $g'_{u'_0}=g_{u'(t'_0)}:(S_{u'_0}, {\bf x}'(u'))\rightarrow (M,{\bf H}'_{f'_0})$.   Since adding-dropping making points does not affect the smoothness, we may assume that
 the number od marked points of $g_u$ and $g'_{u'}$ are the same.
 
 Now fix a    dropping marking map $r=r_k^m$ that selects  $m$   elements ${\bf x}^r(u)$ from the $k$ marked points ${\bf x}(u)$   satisfies the condition that  each free component of an (hence any  ) element $g_{u_0}\in U^{T_1}(g_{u_0})$ is minimally stabilized. This gives a "new" marking, an identification $(S_u)_v\simeq S^2_v=S^2$ of a free component $(S_u)_v, v\in T_1$. Choose a "compatible marking" $r'$ for an (hence any  ) element $g'_{u'_0}\in U'^{T_1}(g'_{u'_0})$ accordingly.

 Let $\Gamma=\Pi_{v\in T_1} \Gamma_v$ acting on the free components of $S_u$.  Here each $\Gamma_v$ is a subgroup of $PSL(2, {\bf C})$  depending on 
 the number of doubles points on the component $S^2_v$. In particular, if $S^2_v$ is stable $\Gamma_v=\{e\}.$
 Now the assumption $[g'_{u'_0}]=[g_{u_0}]$ of the centers above implies that there is a
 ${\widetilde \gamma_0}:S'_{u'_0}\rightarrow  S_{u_0}$ such that $g_{u_0}\circ {\widetilde \gamma_0}=g'_{u'_0}$. Using the identifications (markings) of the free components with (a collection of ) $S^2$(s ) given by $r$ and $r'$, and denoted by 
 $\psi^r$ and $\psi^{r'} $, the map ${\widetilde \gamma_0}=\psi^r\circ \gamma_0\circ (\psi^{r'})^{-1}$ for an element 
 ${\widetilde \gamma_0}\in \Gamma$, unique upto the finite isotropies of $g_{u_0}$.
 Thus we can define the $\Gamma$-action on $U^{T_1}(g_{u_0})$ by a similar formula,
  $\gamma\cdot g_{u}=:g_{u}\circ (\psi^r\circ \gamma\circ (\psi^{r'})^{-1})$. Here
  $\psi^r=\{\psi^r_v, v\in T_1\}$ with $\psi^r_v:(S_u)_v\rightarrow S^2$, similarly for
 $\psi^{r'}.$ Note that for all elements in  $ g_{u}\in U^{T_1}(g_{u_0})$ the domains
  are all the "same"; the parameter $u$ in te notation $\Sigma_u$ or $S_u$ only describes the locations of the distinguished points.  In term of this "action" of $\Gamma$, we have $g'_{u'_0}=\gamma\cdot g_{u_0}$ between the two centers. 
 
  Denote $(\psi^r\circ \gamma\circ (\psi^{r'})^{-1}:S_{u'}\rightarrow S_u$ by ${\widetilde \gamma}$. Let ${\bf x}'^{\gamma_0}(u'_0):= {\widetilde \gamma}_0^{-1}({\bf x}(u)\in S_{u'}$ be the "new" marked points. Note that $g'_{u'_0}({\bf x}'^{\gamma_0}(u'_0))=g_{u_0}\circ {\widetilde \gamma}_0({\widetilde \gamma}_0^{-1}({\bf x}(u))=g_{u_0}({\bf x}(u))\in {\bf H}_{f_0}$.  Now consider  the collection of the $L_k^p$ maps $g'_{u'}$ of type $T_1$  near  $g'_{u'_0}$ with the  constrains on the new marked points: $g'_{u'}({\bf x}'^{\gamma_0}(u'))\in {\bf H}_{f_0}, $
 denoted by  
  $ U'^{T_1}(g'_{u_0}, {\bf H}^{\gamma_0}_{f_0})$.
 
   
   \begin{theorem}
   Let ${\xi}'^{\Phi, \gamma_0} $ be the section ${\xi}^{\Phi} $ viewed in $ U'^{T_1}(g'_{u_0}, {\bf H}^{\gamma_0}_{f_0})$. Then ${\xi}'^{\Phi, \gamma_0} $ is of class $C^{m_0}$.
   \end{theorem}
   
   \proof
   
   Let $G_e^3=\Gamma\times W^{T_1}(f_{t_0}, {\bf H}_{f_0})(=\Gamma\times G_e^1=\Gamma\times G_e^2$ with action of $G_e^2$ first then composing wit the action of $\Gamma$. Then the discussion above shows that ${\xi}'^{\Phi, \gamma_0} $ is the restriction  to the slice $ U'^{T_1}(g'_{u_0}, {\bf H}^{\gamma_0}_{f_0})$ of the $G_e^3$-extension of $\xi^{t_0}$.  Hence it is of class $C^{m_0}$.
   
   \QED 
   
  Applying implicit function theorem in a similar  way  to the proof that $\Psi$ induced by dropping makings is a differomorphism implies the next corollary.
  
   \begin{cor}
   	Let ${\xi}'^{\Phi} $ be the section ${\xi}^{\Phi} $ viewed in $ U'^{T_1}(g'_{u_0}, {\bf H}'_{f'_0})$. Then ${\xi}'^{\Phi} $ is of class $C^{m_0}$.
   \end{cor}
   
   This proves the half of the main theorem  below.

    \begin{theorem}
      Let ${\xi}^{K}:{\widetilde W}(f_{K})\rightarrow {\cal L}^{K}$ be a smooth section satisfying the condition $C_1$ and $C_2$. Then the smooth section ${\xi}$ on $W^{\nu(a_{T_1}), T_1}(g_{t_0}, {\bf H}_{f_0})$ 
      is  of class $C^{m_0}$ viewed in any other such local slices $W^{\nu(a_{T_1}), T_1}(g'_{t'_0}, {\bf H}_{f'_0})$
       or  $U^{ T_1}(g'_{u'_0}, {\bf H}_{f'_0})$ with respect their own smooth structures. Here $[g'_{t'_0}]=[g_{t_0}]$  as unparametrized maps.

       \end{theorem}

   \proof
    
    The other half essentially  follows from  the obvious "extension  " of  the action of $G_e^3$ defined above by composing further  the identifications given by inverse of 
   $\phi'_{t'}$ (with $t'$ and $t$ corresponding to each other) first, then the inverse of $\lambda^{t_0'}_{t'}$. Indeed, let ${ {\xi}}'$ be the section $\xi$ viewed in
    $W^{T_1}(f_{t_0}, {\bf H}_0)$. 
Then upto the effect of  using different (but fixed ) markings, ${ {\xi}}'$ is the restriction to $W^{T_1}(f_{t_0}, {\bf H}_0)$ of the new $G_e^3$-extension
 of $\xi^{t_0}$.   A similar argument to the proof of the corollary above will eliminate the effect of  different makings so that ${ {\xi}}'$ is of class $C^{m_0}.$
   
  \QED

\end{document}